\documentclass{amsart}
\usepackage{amssymb,latexsym, amsmath}
\newtheorem{theorem}{Theorem}
\newtheorem{lemma}{Lemma}
\newtheorem{definition}{Definition}

\newtheorem{proposition}{Proposition}
\newtheorem{corollary}{Corollary}

\bibliography{samparatb}
\begin{document}

\title{$ \:\:\:\:\:\:$Isometric Equivalence of Isometries on $ H^p $}

\maketitle

\textbf{$\:\:\:\:\:\:\:\:
\textbf{Joseph \:\:A.\:\: Cima} \:\:\:\:\: and \:\:\:\:\: \textbf{Warren \:\:\:R.\:\:  Wogen \:\:\:}
\: $} \vskip 0.4in

Abstract.We consider a natural notion of equivalence for bounded linear operators on $H^p,$ for $p\neq 2.$ We determine which isometries of finite codimension are equivalent. For these isometries , we classify those which have the Crownover property.

  \vskip 0.2in

\section{Introduction}

  If $ A$ and $B$ are bounded linear operators on a Hilbert space, then $A$ and $B$ are unitarily equivalent if $ B= U A U^{*}$ for
 some unitary operator $U.$ One can then view $A$ and $B$ as abstractly the same operator.
 In the general Banach space setting one can replace unitary equivalence by using onto isometries of the Banach space
 being considered. So if $X$ is a Banach space and $ \mathbb{B}(X)$ are the bounded linear operators on X, then we will
 say that  $A$ and $B$  in $ \mathbb{B}(X)$   are isometrically equivalent if $B = U A U^{-1}$ for some onto isometry
 $U$ in $ \mathbb{B}(X).$  In this case we write $ A \thickapprox B.$ The utility of this notation will of course depend
 on specific properties of the space $X$ and its onto isometries. The Banach spaces considered in this note are the classical Hardy spaces $H^p,$ for $p \neq 2.$
The onto isometries of $H^p$ have been determined (see Theorem A) , and we will classify some familiar operators
on $H^p$ up to isometric equivalence. This work is  motivated by some questions of J. Jamison. In particular
he asked which isometries on $H^p$ are equivalent to the shift ( see Corollary $1$).  \vskip 0.2in

\section{Preliminaries}

  In this paper we consider the Banach spaces $ H^p$ of the unit disc D , for $\:\:\:\:\:\:$
  $1\leq p < \infty, p \neq 2.$
 Recall that $H^p$ consists of the analytic functions $f$ on D for which

 \vskip 0.1in
  $$ \:\:\:\:\:\:\: Sup_{0<r<<1} \int_T |f(r\zeta)|^p dm(\zeta)
$$is finite , $dm$ the usual Lebesgue measure on the unit circle T (See Duren{$2$]).
 \begin{definition}
 $\mathbb{I}(H^p)$ will denote the onto isometries of  $H^p\:, \: p \neq 2 . $
 \end{definition}
 \vskip 0.1in
 $\mathbb{I}(H^p)$ is a group under the usual operator multiplication. The description of   $\mathbb{I}(H^p)$  for $p=1$ is due to de Leew,et.al. [$4$] , and for $1<p <\infty\:\: ,\:  p \neq 2\:,$ to Forelli [$6$]  and is given in Theorem A below.
  \vskip 0.1in
  $2010$ \textit{Subject Classification} Primary $47$B$32, 33$, $30$J$05.$\vskip 0.1in
  \textit{Key words and phrases.} Hardy spaces, Isometries.
 \begin{definition}
 Let $\mathbb{A}$ be the collection of holomorphic automorphisms of the unit disc D.
 That is,
  $$
  \mathbb{A} = \{ \phi(z) =  \frac{ \lambda(z-a)}{1 - \overline{a} z};  \:\:\: a \in D, \:\:\: |\lambda|=1\}
  $$
\end{definition}

\vskip 0.2in

 $\mathbb{A}$ is a group under composition with identity e, where $e(z) = z.$ Following [$8$], we write $\phi_n$ for the n-fold composition of $\phi$ with itself. In addition, we denote the (compositional) inverse of
 $\phi$ by $\phi_{-1}.$  Note that $(\phi_n)_{-1},$ the inverse of $\phi_n,$ is just $(\phi_{-1})_n,$ the n-fold
 composition of $\phi_{-1} $ with itself, which we denote by $\phi_{-n}.$

 \begin{definition}
 For $\phi,$ and $\psi$ in $\mathbb{A}$ we say that $\phi$ is conjugate to $\psi$ if there is an $\eta \in \mathbb{A}$ with
  $\phi = \eta \circ \psi \circ \eta_{-1}. $ The conjugacy class of $\phi$ is denoted  $\mathbb{C}(\phi)$ and so
  $\mathbb{C}(\phi) = \{\eta \circ \phi \circ \eta_{-1} \: \:: \: \eta \in \mathbb{A}.\}$
  \end{definition}

  \vskip 0.3in

 The following proposition is  well known and we include it for notational purposes.

  \begin{proposition}
 Suppose $(b_k)$ is a sequence of functions in $\mathbb{A}$ and that $b_k(a_k)=0\:\: \forall \:\: k \geq 1.$

\vskip 0.2in
 (i) If $\prod_{k=1}^{\infty} b_k(z) $ is a convergent Blaschke product, then $(a_k)$ is a Blaschke sequence. That is,
 $\sum_{k=1}^{\infty} ( 1 - |a_k|)  < \infty.$ \vskip 0.2in
 (ii) Conversely, if $(a_k)$ is a Blaschke sequence then there exists $(\lambda_k)_{k=1}^{\infty} \in T$ so that
 $\prod_{k=1}^{\infty} \lambda_k b_k(z)$   converges.
\vskip 0.3in
 \end{proposition}

 \section{Isometric Equivalence}

   A theorem of Forelli [$6$ ] describes all isometries of $H^p$ onto $H^p.$
  For $\phi \in \mathbb{A}$ and $d\in T$ Forelli showed that  the map
  $$
   f \mapsto d (\phi'(z))^{1/p} f\circ\phi\:\:\: $$
    is in $ \mathbb{I}(H^p) ,$
   and that all onto isometries have this form. Theorem A below is a  restatement of this result.

  \vskip 0.2in

   If
   $$\phi(z) =  \frac{\lambda(a-z)}{1 - \overline{a} z},$$
   \vskip 0.2in
   then
   $$ \phi'(z) =  \frac{\lambda(1 - |a|^2)}{ (1 - \overline{a} z)^2},
   $$
   so that the choice of a branch of the pth root function that will make $(\phi'(z))^{1/p}$ analytic will
   depend on $\lambda.$ It is useful to set our notation so that we always use the principal branch given
   by $(r \exp(\imath \theta))^{1/p} =r^{1/p} \exp(\imath \theta/p)\:\: , \:\: - \pi < \theta < \pi\:, \: r>0.$
   Now $(\overline{\lambda} \phi'(z)) = \frac{(1 - |a|^2)}{ (1 - \overline{a} z)^2} $ has positive real part
   so $(\overline{\lambda} \phi'(z))^{1/p}$ is analytic on D.

  Let $ C_{\phi} f = f \circ \phi$ denote composition by $\phi$ and for $F\in H^{\infty}, M_F $ denotes
   multiplication by F. Finally, let
   $$
    U_{\phi} = M_{((\overline{\lambda})\phi')^{1/p}} C_{\phi}.
    $$
    \vskip 0.2in
  Forelli's result can be stated as
 \vskip 0.2in
\textbf{Theorem} $A.$  $\mathbb{I}(H^p) = \{\rho U_{\phi}\: : \phi \in \mathbb{A} \:, |\rho|=1 \}$

\vskip 0.2in
   From this result we see that $\mathbb{I}(H^p) $ is determined by $\mathbb{A}.$ We will examine the
  relation between the group structure of $\mathbb{A}$ and $\mathbb{I}(H^p).$
\vskip 0.2in

   We note that for $p=2,$  the  operators of the form $\rho U_{\phi}$ in Theorem $A$ are of course unitary
  operators on $H^2$ which are tied to the analytic structure of $H^2.$ These unitaries are a small subgroup of the
  full unitary group on $H^2.$

  \begin{lemma}
  Let $ \phi ,$ and $ \psi \in \mathbb{A}.$ Then

  $a) U_{\phi} U_{\psi} = \rho U_{\psi\circ\phi} $  for some $\rho \in T,$ which depends on $\phi$ and $\psi.$
  \end{lemma}

 $b) U_{\phi}^{-1} = U_{\phi_{-1}}.$
 \vskip 0.1in
 Proof. Suppose that $\phi(z) = \frac{\lambda_1 (z-a_1)}{1 - \overline{a_1} z} $ and
  $\psi(z) =  \frac{\lambda_2(z - a_2)}{1 - \overline{a_2} z}.$ Then
  $$
  \psi \circ \phi(z) =  \frac{\lambda_3(z - a_3)}{1 - \overline{a_3} z}
  $$
 for some $\lambda_3 \in T, a_3 \in D.$ Note that $C_{\psi} C_{\phi} = C_{\phi\circ\psi}.$ Also if $F\in H^{\infty},$
 then $C_{\phi} M_F = M_{F\circ \phi} C_{\phi}.$ Thus
 $$
 U_{\phi} U_{\psi} = M_{(\lambda_1\:\: \phi')^{1/p}} C_{\phi} M_{(\lambda_2\:\:  \psi')^{1/p}} C_{\psi \circ \phi}=
 (\overline{\lambda_1} \phi')^{1/p} (\overline{\lambda_2} \psi'\circ \phi)^{1/p} C_{\psi \circ \phi}.
 $$
 \vskip 0.2in
 But
 $$
  U_{\psi \circ \phi} = ( \overline{\lambda_3} (\psi\circ\phi)' )^{1/p}C_{\psi\circ\phi} =
  (\overline{\lambda_3}(\psi'\circ\phi) \phi')^{1/p} C_{\psi\circ\phi}
  $$
 \vskip 0.2in
  So one sees that $U_{\psi\circ\phi}$ is a unimodular multiple of $U_{\phi} U_{\psi},$ and a) is proven.  \vskip 0.2in
  For part b) we recall that
  $$
   \phi_{-1}(z) =  \frac{\overline{\lambda_1}(z + \lambda_1 a_1)}{(1 + \overline{\lambda_1 a_1}z)}.
   $$
 Take $\psi = \phi_{-1}$ in the last proof. Thus
 $$
 U_{\phi} U_{\phi_{-1}}= (\overline{\lambda_1} \:\phi')^{1/p}(\lambda_1 \phi_{-1}'\circ \phi)^{1/p} C_{\phi_{-1} \circ \phi}=I\:\:\:\:\:\:\#
 $$
\vskip 0.2in

  Remark. The value of the constant $\rho$ in Lemma $1 $ a)  will not be needed in our work, but it can of course be
explicitly   computed.   For $\phi $ and $\psi$ as in the proof of Lemma $1,$ one can show that $\rho = \exp i \theta,$
where $\theta = arg (1 + \overline{\lambda_1 a_1} a_2)^{2/p}.$

\vskip 0.2in
  We now describe all $ S \in \mathbb{I}(H^p)$ which are isometrically equivalent to a fixed $U_{\phi} \in \mathbb{I}(H^p).$\vskip 0.2in
\begin{proposition}
 $\:S \thickapprox U_{\phi} \Leftrightarrow $ there exists $\eta \in \mathbb{A} $ and $\rho \in T$ so that
 $$
 S = \rho U_{\eta \circ \phi \circ \eta_{-1}}.
 $$
\end{proposition}
\vskip 0.1in
Proof: $ \: S \approx U_{\phi}  \Leftrightarrow $ there exists $\eta \in \mathbb{A}$ so that
$$
U_{\eta_{-1}} U_{\phi} U_{\eta} = S.
$$
 But
 $$
 U_{\eta_{-1}} U_{\phi} U_{\eta} = U_{\eta_{-1}}(\rho_1 U_{\eta \circ \phi})= \rho_2 \rho_1 U_{\eta \circ \phi \circ \eta_{-1}}
 $$
 where $\rho_1, \rho_2 \in T$ as in Lemma $1 a.\:\:\:\:\:\:\:\:$         \#   \vskip 0.1in
 So Proposition $2$ states that
 $$
 \widetilde{\psi} \in \mathbb{C}(\phi) \Leftrightarrow U_{\phi} = \rho U_{\widetilde{\psi}}
 $$
 for some $\rho \in T.$

 \vskip 0.3in
 We focus on the isometries of $H^p$ into $H^p.$ The most familiar example is the shift, $M_z,$ on $H^p.$ The range of $M_z$ is $z H^p$ so is of codimension one.   For example, which $ \:S \in \mathbb{B}(H^2)$ satisfy $\:S \thickapprox M_z?$ We will in fact classify
 all  finite codimension  isometries up to isometric equivalence.  We give the details of our results for the case for codimension one isometries and the codimension n case follows similarly.

 \vskip 0.2in
   Note that $M_z$  has the additional property that
 $$
  \bigcap_{n=1}^{\infty}(M_z)^n H^p = (0).
  $$

  \vskip 0.2in
  \begin{definition}
   A codimension one isometry  $S$ on $H^p$ is called Crownover (see [$2$] , [$7$]) ,  if $
 \bigcap_{n=1}^{\infty} S^n H^p = (0).$
  \end{definition}

  We will also classify such isometries up to isometric equivalence. \vskip 0.2in

  \vskip 0.3in

  \section{Finite Codimensional Isometries }

  We will now state Forelli's theorem [$6$, Theorem $1$] describing all isometries of $H^p,\: p \neq 2.$
  \vskip 0.2in

  \textbf{Theorem} $B.$ $S$  is an isometry of $H^p,  \: 1 \leq p < \infty\:, \: p \neq 2$ iff  $S = F f(\phi) $
  for some $\phi$ inner and an $F \in H^p$ which is related to $\phi.$

\vskip 0.2in
   The precise relationship between $F$ and $\phi$ can be found in [$6$] and is not needed in the work that follows. We will provide a simpler description of the isometries of finite codimension.

  \vskip 0.2in
 \begin{lemma} Suppose that  $ T = M_F C_{\phi}$  as in $Theorem\:\: B$ and that the inner function $\phi$ is  not in   $\mathbb{A}.$ Then the T has infinite codimension.
 \end{lemma}

 Proof: We will modify the  proof of [$1$ , Lemma $3.6$] ] . Let $K_b(z)= \frac{1}{1-\overline{b} z}. $

 Since $\phi$ is an open map which is not univalent, we can choose sequences $(a_n),(b_n)$ in D so that $ \phi(a_n)= \phi(b_n) = c_n\:\:\forall\: n.$ F is not the zero function
 so we can also assusme $F(a_n)\neq 0$ and $F(b_n)\neq 0 \:\:\forall\: n.$  Let $g_n = \overline{F(b_n)} K_{a_n} - \overline{F(a_n)} K_{b_n}. $
 Since the kernels are linearly independent the functions $g_n$ are linearly independent. The $g_n $ are in $H^{\infty}$ and so induce linear functionals $\Lambda_n$ on $H^p$ which are linearly independent and  satisfy
 $$
 \Lambda_n(Tf)= \Lambda_n (F f\circ \phi) = F(b_n)F(a_n)c_n - F(b_n)F(a_n)c_n =0,
 $$
 for all $f \in H^p.$ Hence,
 $$
 \bigcap Ker (\Lambda_j) \supset T(H^p).
  $$

  Thus $\{g_n\}$ is a linearly independent set whose span intersects $ T H^p$ only at $(0).$  This implies  T has infinite codimension. $\:\:\:\:\: \#\:\:$

 \vskip 0.3in
   Hence  we need only consider isometries of the form $M_F C_{\phi}\:\:,$ where $\phi \in \mathbb{A}$ and $F \in H^p.$
  If $\phi(z) = \frac{\lambda(z-c)}{1 - \overline{c} z} ,$ then $  M_F C_{\phi}\: = M_{\frac{F}{(\overline{\lambda} \: \phi')^{1/p}}} \:\: U_{\phi} .$
  Since $U_{\phi} \in \mathbb{I}(H^p) , $  it follows that $M_{\frac{F}{(\overline{\lambda} \: \phi')^{1/p}}} \:\: $
  must be isometric. This means that $M_{\frac{F}{(\overline{\lambda} \: \phi')^{1/p}}} \:\: $ is an inner function,
  which we label as $\Psi.$
   \vskip 0.2in

    Clearly $ M_{\Psi} U_{\phi}$ has the codimension of $M_{\Psi}.$ The codimension  is $ n < \infty \Leftrightarrow \Psi$
  is an n-fold Blaschke product. In this case we write $\Psi \in \mathbb{A}_n.$  In particular $\mathbb{A}_1 = \mathbb{A}.$

    We have shown that the set of isometries of codimension n is given by
    $$
     \mathbb{I}_n(H^p) = \{ M_{\Psi} U_{\phi} \:\: : \: \phi \in \mathbb{A}\:, \Psi \:\in \mathbb{A}_n \:\}.
     $$
  In most of what follows, we focus on the isometries
   $$
    \mathbb{I}_1(H^p) = \{ M_{\psi} U_{\phi} \: : \phi, \psi \in \mathbb{A}\}
    $$
   of codimension one.

\vskip 0.3in
\begin{theorem}
Let $S_1 = M_{\psi} U_{\phi} \in \mathbb{I}_1(H^p).$  If $S_2 \in \mathbb{I}_1(H^p),$ then $S_2 \approx S_1 \Leftrightarrow
\exists \eta \in \mathbb{A} $ and $ \rho \in T$ so that $S_2 = M_{\rho \psi \circ \eta} U_{\eta_{-1} \circ \phi \circ \eta}.$
\end{theorem}
 \vskip 0.1in
 Proof: $S_2 \thickapprox S_1 \Leftrightarrow \exists \eta \in \mathbb{A} $ so that $U_{\eta_{-1}} S_1 U_{\eta} = S_2.$ But
 $$
 U_{\eta} S_1 U_{\eta{-1}} = U_{\eta} M_{\psi} U_{\phi}U_{\eta_{-1}}
 = M_{\psi \circ \eta} U_{\eta} (\rho_1 U_{\eta_{-1}\circ \phi})
 = M_{\psi \circ \eta} \rho_1 (U_{\eta} U_{\eta_{-1} \circ \phi})$$
  $$=
  M_{\psi \circ \eta} \rho_1 \rho_2 U_{\eta_{-1} \circ \phi \circ \eta}=
  \rho  M_{ \psi \circ \eta} U_{\eta_{-1} \circ \phi \circ \eta}.
 $$
 Here $\rho_1$ and $\rho_2$ are the unimodular constants that arise in Lemma $1$ a, and $\rho = \rho_1 \rho_2.\:\:\:\:\:\:$       \#

 \vskip 0.2in
  With $e(z)=z$ note that if $\psi \in \mathbb{A},$ then $M_{\psi} = M_{\psi}U_{e}$ has codimension one .
  \begin{corollary}
   If $\psi \in \mathbb{A}$ and $S \in \mathbb{I}_1(H^p),$ then $\:S =  M_{\widetilde{\psi}}$ for some $\widetilde{\psi} \in \mathbb{A}.$
  \end{corollary}

Proof: $ S \thickapprox M_{\phi} \Leftrightarrow \exists \eta \in \mathbb{A}$ so that
$$
 \:S= U_{\eta} M_{\psi} U_{\eta_{-1}} = M_{\psi \circ \eta} U_{\eta} U_{\eta_{-1}}= M_{\psi \circ \eta}.
$$
Finally, note that $\{ \psi \circ \eta \:\: : \:\: \eta \in \mathbb{A}\} = \mathbb{A}
 \:\:\:\:\:\:\: \#$
\vskip 0.3in

 We remark that the above result shows that $\: S \thickapprox M_z \Leftrightarrow S = M_{\psi}$ for some $\psi \in
 \mathbb{A}.$ We now generalize the last corollary. Fix $\phi \in \mathbb{A} $ and consider when $M_{\psi} U_{\phi} \thickapprox  M_{\widetilde{\psi}} U_{\phi}.$ Corollary $1$ settles the question if $\phi = e.$

 \vskip 0.2in
  So suppose $\eta \in \mathbb{A}$ and that $U_{\eta} (M_{\psi} U_{\phi})U_{\eta_{-1}} = M_{\widetilde{\psi}} U_{\phi}.$
  The left side simplifies to
  $$
  M_{\psi \circ \eta} U _{\eta} U_{\phi}U_{\eta_{-1}} = \rho M_{\psi \circ \eta} U_{\eta_{-1} \circ \phi \circ \eta},
   $$
and with the notation  $\widetilde{\psi} = \rho \psi\circ \eta$ and $\phi = \eta_{-1}\circ \phi \circ \eta $ we have our equality
$$
U_{\eta} (M_{\psi}U_{\phi})U_{\eta_{-1}} = M_{\widetilde{\psi}}U_{\phi}.$$ It follows that $\phi \circ \eta = \eta \circ \phi,$ so $\phi$ and $\eta$ commute. Thus we have

\vskip 0.1in
\begin{corollary}
 $M_{\psi} U_{\phi} \thickapprox M_{\widetilde{\psi}}U_{\phi} \Leftrightarrow \widetilde{\psi} = \rho \psi\circ \eta$ for some $\eta \in \mathbb{A}$ with $\eta$ commuting with $\phi$ and $\rho \in T$ satisfying
 $$
 U_{\eta} U_{\phi} U_{\eta_{-1}} = \rho U_{\eta_{-1} \circ \phi \circ \eta}.$$

\end{corollary}

\vskip 0.1in
Remark: We will discuss in Section $7$  the classification of the automorphisms commuting with a fixed $\phi \in \mathbb{A}.$
\vskip 0.2in

\vskip 0.2in
 Recall that $ \forall \:\: \psi \in \mathbb{A} , M_{\psi} $ is   a Crownover shift. That is
$$
\bigcap_{n=1}^{\infty} (M_{\psi})^n H^p = \bigcap_{n=1}^{\infty} (M_{(\psi)^n)}) H^p = (0).
$$
Given a $ S = M_{\psi} U_{\phi} \in \mathbb{I}_1(H^p)$ when is $S$ Crownover?
Now $$
 \:S^2 = (M_{\psi} U_{\phi})((M_{\psi}U_{\phi}) = M_{\psi} M_{\psi \circ \phi} U_{\phi} U_{\phi}= \rho M_{\psi} M_{\psi \circ \phi} U_{\phi_2},$$ for some $\rho \in T.$  Iterating we have
$$
S^n = (M_{\psi} U_{\phi})^n = \rho  M_{\psi}M_{\psi\circ\phi}\centerdot\centerdot\centerdot\centerdot M_{\psi \circ \phi_{n-1}}U_{\phi_{n}},
$$
where $\rho \in T$ depends on n.

Now $U_{\phi_n}$ is onto , so $S^n H^p = B_n H^p ,$ where $B_n$ is the Blaschke product $\prod_{k=0}^{n-1} b_k,$ where
$$b_k = \psi \circ \phi_k.
$$

  Note that $b_k$ is merely the kth term of the sequence $(\psi\circ\phi_k)$ and does not represent the kth iterate of $b.$
\vskip 0.3in
 It follows that $\bigcap_{n=1}^{\infty} S^n H^p = \bigcap_{n=1}^{\infty}B_n H^p .$ If this intersection  contains an $f \neq 0$ then each $b_k$ is a factor of $f$ so by Proposition $4$ there is a Blaschke product of the form $B=\prod_{k=0}^{\infty} \lambda_k b_k$ such that $ \bigcap_{n=1}^{\infty} S^n H^p = B H^p.$ Thus the zeros of $(b_k)_{k=0}^{\infty} $ form a Blaschke sequence. The above discussion shows that

 \begin{theorem}
  $M_{\phi}U_{\phi}$ is Crownover $ \Leftrightarrow $ the sequence of zeros of $(\psi \circ \phi_k)_{k=0}^{\infty}$
  is not a Blaschke sequence.
  \end{theorem}

We will elaborate on this result in the  next   section using the fixed point structure of $\phi.$

At this time we maintain the terminology as above , assuming that $S = M_{\psi} U_{\phi}$ and that $B= \prod_{k=0}^{\infty} \lambda_k b_k$ is an infinite Blaschke product, with $\lambda_0= 1.$ Note that $BH^p$ is an invariant subspace for $S.$ We will show that $S|_{B H^p} \in \mathbb{I} (B H^P)$. First note that
$$
M_{\psi} C_{\phi}B= M_{\psi} C_{\phi} \prod_0^{\infty} \lambda_k b_k = \psi  \prod_0^{\infty} \lambda_k b_{k+1}
$$
So if $ g \in H^p,$ then
$$
 S Bg = M_{\psi}U_{\phi} B g= B U_{\phi} g,
$$
and $S|_{B H^p} $ is onto $B H^p.$

Lastly, we note that $S|_{B H^P} $ is isometrically equivalent to $U_{\phi}.$  Let $ V \: :
 \:\: H^p \rightarrow   B H^p$ be the isometry defined by

$$
 V _g = B g\:\: , \:\: g \in H^p.
 $$
 Then
 $$
  g \in H^p \Rightarrow(S|_{B H^p}) V_g = S(B g) = B U_{\phi} g,
  $$
  so that $S|_{B H^p} \approx U_{\phi}.$

\vskip 0.3in

Remark: We now consider the  case as above  but with $p=2.$  The Wold decomposition for the isometry $ S$ (see [$9,$ Th.$ 1.1 $]),  is easy to exhibit. Namely,
$ H^2 = B H^2 \bigoplus (B H^2)^{\bot}$ is  a direct sum of invariant subspaces  of $S.$ $S|_{B H^2}$ is unitary and is in fact unitarily equivalent to $U_{\phi},$ while $S|_{(B H^2)^{\bot}}$ is a unilateral shift. If $\psi(z)  = \frac{\mu (z-b)}{1 - \overline{b} z} ,$ then span of the kernel $K_b$  is a wandering subspace for the shift.

\vskip 0.3in

\section{The Crownover Property}

\vskip 0.2in
 Each $\phi \in \mathbb{A}, \phi \neq  e , $ can be classified as elliptic, hyperbolic, or parabolic according to its fixed points in $\overline{D}.$ See [$1$] or [$8$] for more detail.

  \begin{definition}
    $\phi \in \mathbb{A}\:,\: \phi \neq e \:,$ is elliptic if $\phi$ has a fixed point, say a, in D. Let $\mathbb{E}(a) =\{\psi \in \mathbb{A} \:\: : \:\: \psi(a) = a , \:\: \psi \neq  e \}\:$ denote the set of all elliptic automorphisms of D that fix a.

  \end{definition}

  Choose $\eta \in \mathbb{A} $ with $\eta(a) = 0$ and note that $\eta \circ \mathbb{E}(a)\circ \eta_{-1}$ is the set of
  nontrivial  rotations of D.

 \vskip 0.2in
\begin{definition}$\phi$ is parabolic if it has only one fixed point , say w. In this case $w\in T$ and $\phi'(w)=1.$ Of course w is also
 the unique fixed point of $\phi_{-1}.$ w is attractive for $\phi$ (and for $\phi_{-1}).$ That is, for all $c \in D, \phi_n(c) \rightarrow w$ and $ \phi_{-n}(c) \rightarrow w.$
\end{definition}

\vskip 0.2in
\begin{definition}
$\phi \in \mathbb{A}\:\:, \: \: \phi \neq e $ is called hyperbolic if $\phi$ has two distinct fixed points, say $w_{1},$ and $w_{2},$ on T. In this case one of the fixed points , say $w_{1},$ is the attractive  fixed point for $\phi.$ Also $w_{2}$ is the attractive fixed point for $\phi_{-1}.$  Further, $\phi'(w_1) <1 $  and $\phi'(w_2) >1.$
\end{definition}

\vskip 0.2in
 Let
 $$
  \mathbb{H}(w_{1}\:,\:w_{2})= \{\psi \in \mathbb{A}\:, \phi  \not = e \::\psi (w_{1}) = w_{1} \: , \: \psi(w_{2}) = w_{2} \}
  $$
  be the collection of hyperbolic automorphisms that fix $w_1 $ and $w_2.$  If $w_1^{'}, w_2^{'}$ is another pair of distinct points on T and $\eta \in \mathbb{A}$ is chosen so that $\eta(w_1)=w_1^{'}, \eta(w_2)=w_2^{'},$  then
 $\mathbb{H}(w_1^{'},w_2^{'})= \eta \circ \mathbb{H}((w_1,w_2)\circ \eta_{-1}.$  Thus all of these sets are conjugate.
 \vskip 0.2in
 As an example , take $w_1=-1\:\:,\:\: w_2 = 1.$ Then one can show $\mathbb{H}(-1,1) = \{ \psi_r \:\:; \:\: -1 < r <0 \:or \: 0 < r <1\},$
 where $\psi_r (z) = \frac{z-r}{1- r z}.$

 \begin{proposition}
 If $\phi $ is elliptic and $ \psi \in \mathbb{A},$  then $M_{\psi} U_{\phi}$ is Crownover.
 \end{proposition}

 Proof: Since $\phi$ is conjugate to a rotation , it is routine to check that the zeros of  $\psi \circ \phi_n$ lie on a circle in D and hence can not be a Blaschke sequence. $\:\:\:\;\:\:\:\#$

 \vskip 0.4in

\begin{proposition}
If $\psi \in \mathbb{A}$ and $\phi$ is hyperbolic , than $M_{\psi}U_{\phi}$ is not Crownover.
\end{proposition}

\vskip 0.2in
Proof: It is easy to check that if $\phi$ is hyperbolic and $c \in D$ that $\sum(1 - |\phi_n(c)|) < \infty.$  (See [$8,p 85,\#6.$) Suppose $\psi \circ \phi_n(a_n) = 0 \:\:\forall  \:\:n \geq 0.$ Then $\phi_n(a_n) = \psi_{-1} (0),$ and $a_n= \phi_{-n}\circ \psi_{-1}(0).$  But $ \phi_{-1}$ is also hyperbolic , so $\sum(1 - | \phi_{-1}(\psi_{-1}(0)|)= \sum(1 - |a_n|)< \infty $ and $(a_n)$ is a Blaschke sequence.  $\:\:\:\: \: \# \:\:\:$

 \vskip 0.2in

 Our goal is to show that if $\phi,\psi \in \mathbb{A}$ with $\phi$ parabolic, then $M_{\psi} U_{\phi}$ is not Crownover.
 That is , the zeros of $(\psi \circ \phi_n)$ form a Blaschke sequence, just as in the case that $\phi$ is hyperbolic.
 \vskip 0.2in
 \begin{definition} For $w \in T,$ let  $\mathbb{P}(w)$ be the collection of parabolic automorphisms that fix w.
\end{definition}
 It is easy to see that the sets  $\mathbb{P}(w),$  $ w \in T \:,$  are conjugate to each other.
 So we first consider $\mathbb{P}(1).$

 A computation will show that
$$
 \phi(z) =\frac{\lambda(z- a)}{1- \overline{a} z} \in \mathbb{P}(1) \Leftrightarrow  \phi(1)=1=\phi'(1).
 $$

Solving for a and $\lambda,$ we see that
$$
 |a - 1/2| = 1/2 \:\:,\:\: a \neq 0,1$$
 and  that $\lambda = \frac{1- \overline{a}}{1 - a}.$  So $\: a - 1/2  = (c/2) ,$ where $ c \in T,c \neq \pm 1$ and thus
 $\lambda = \frac{1 - \overline{c}}{1 - c} = \frac{-1}{c}.$
 Using these equalities we can write $\phi$ in the form
 $$
 \phi(z) = \frac{ 1 + c -2z}{2c -(1 + c) z},
 $$
 which we write as $\phi_c(z).$

\vskip 0.2in
So we have
 $$
 \mathbb{P}(1)= \{ \phi_c(z) = \frac{  1 + c -2z}{2c -(1 + c) z}\:\: : \:\: c \in T\:,c \neq \pm1\}\:.\:
 $$

\vskip 0.2in

The functions $\phi_{\imath}$ and $\phi_{-\imath}$ play a special role in what follows.
\vskip 0.2in
  Observe that if $\phi \in \mathbb{P}(1)$ and if $\psi \in \mathbb{A}$ with $\psi(1)=1,$ then
  $$
  \psi \circ \phi \circ \psi_{-1} \in \mathbb{P}(1).
  $$
  Here $\psi$ could be hyperbolic. Our approach is to conjugate $\phi_{\imath}$ (or $\phi_{-\imath}$) by automorphisms $\psi_r \in \mathbb{H}(-1,1),$ discussed after Definition $7.$ We note that the inverse of $\psi_r$ is $\psi_{-r}.$

 \begin{proposition}
  Let $\phi_c\:(z) = \frac{1 + c - 2z}{ 2c -(1 + c)z}$  where $c \in T\: \: ,\: c \neq \pm i.$ If $ \Im (c) >0, $  then $ \exists \:\:r\:\:  \in (-1,1)\:,\: r \neq 0 $ and $ \psi_r \in \mathbb{H}(-1,1)$ so that $\phi_c = \psi_r \circ \phi_i \circ \psi_{-r},$ while if $\Im (c) < 0\: \:,\:\: \exists\:\:$ another $\:\:r \:\:\ \in (-1,1)\: , r \neq 0 $ so that $\phi_c = \psi_r \circ \phi_{-i} \circ \psi_{-r},$
  \end{proposition}

Proof: Let $r \in (-1,1)\::\:r \not = 0.$ Then $\psi_{r} \circ \phi_{\imath} \circ \psi_{-r} \in \mathbb{P}(1),$
so if $z_r$ is  the zero of $ \psi_r \circ \phi_i \circ \psi_{-r},$ then $|z-1/2|=1/2 .$ Letting $c_r = 2 z_r -1,$ it suffices to show that $\{ c_r\::\: -1 <r < 1\}$ is the upper half semicircle of T.

  \vskip 0.2in
  A careful computation shows that
  $$
  \psi_r \circ \phi_i \circ \psi_{-r} (z) = \frac{(1-r)^2 - ((1-r)^2 - i ( 1 - r^2))z}{ (1 - r)^2 + i ((1-r^2) - (1-r)^2)z}
  $$
  so that
  $$
   z_r = \frac{(1-r)^2 + i (1 - r^2)}{2(1 + r^2)}
   $$
   Thus
   $$
    c_r = 2 z_r -1 = \frac{-2r + i (1 - r^2)}{1 + r^2}.
    $$
One checks that as r goes from $-1$ to $1\:,\: c_r$ traces  out the required semicircle. A similar computation for
$\phi_{-1}$  yields a $c_r$ that traces out the lower semicircle of T.     \#

\vskip 0.2in
\begin{theorem}
$\mathbb{C}(\phi_i)\bigcup \mathbb{C}(\phi_{-\imath }) = \mathbb{P} ,$ the collection of all parabolic automorphisms.
\end{theorem}

Proof:  $\mathbb{P} = \bigcup_{w \in T} \: \mathbb{P}(w).$ Given $w \in T,$ choose $\eta \in \mathbb{A} $ so that
$\eta(1) = w.$ Then we have
$$
\eta \circ \mathbb{P}(1) \circ \eta_{-1} = \mathbb{P}(w).
$$
Thus each $\psi \in \mathbb{P}(w) $ is conjugate to some $\phi \in \mathbb{P}(1),$ and our previous result shows that
$\phi $ is conjugate to $\phi_{i}$ or to $\phi_{-i}.$ Thus $\psi$ is conjugate to $\phi_i$ or to $\phi_{-i}.\:\:\:\:\#\:\:$
\vskip 0.3in

We will now examine  the case where $c= \pm i\:\:, \:\phi_i (z)$ and its inverse $\phi_i = \phi_{-i}.$\vskip 0.2in
\begin{lemma}
The zeroes of the  iterates of $\phi_{\imath}$ (and those of $\phi_{-\imath}$) form a Blaschke sequence.
\end{lemma}

  Proof: Multiplying each coefficient of $\phi_i$ by $(1-i)/2$ shows that
$$
 \phi_i(z)= \frac{1 - (1-i) z}{ 1 + i -z}.$$
 An easy computation shows that
 $$
  \phi_i\:\circ \phi_i (z) = \frac{2 - (2-i)z}{(2 + i) - 2z}
  $$
  and by induction we see that the nth iterate  is given by
  $$
  (\phi_i)_n \:(z) = \frac{n- (n-i) z}{n+i -nz}
  .$$
 Thus $a_n = n/(n-i)$ is the zero of $(\phi_i)_n \:\:  \:\: , \:\:\: |a_n|^2 = n^2/(n^2 + 1).\:\:$
 So $\sum (1 - |a_n|^2) < \infty,$ and $(a_n)$ is a Blaschke sequence. Essentially the same argument
 shows that the zeroes of $(\phi_{-i})_n $ also form a Blaschke sequence.$\:\:\:\#\:\: $

 \vskip 0.2in
 \begin{lemma}
 Suppose that $\phi \in \mathbb{A}$ with $\phi_n (a_n) =0 \: \; \forall \: n  \Rightarrow (a_n) $ is a Blaschke sequence. Then

 i) If $\psi \in \mathbb{A}$ and $\psi \circ \phi_n(b_n) = 0\:\:, \: \forall  n, $ then $(b_n)$ is a Blaschke sequence.

 ii)If  $\widetilde{\phi}  \in \mathbb{C}(\phi) $ and if  $ (\widetilde{\phi})_n(c_n) =0 \:, \: \forall n,$ then $(c_n)$ is a Blaschke sequence.
 \end{lemma}

 Proof: For i) assume $(a_n)$ is a Blaschke sequence and consider $ \psi \circ \phi_n(b_n) =0,$
   with $\phi_n(z) = \frac{\lambda_n(z - a_n)}{1 - \overline{a_n} z}.\:\:\:$
  So $b_n = \phi_{-n} \circ \psi_{-1} (0).$ Let $\alpha = \psi_{-1} (0) \in D, $  and $|\phi_{-n} (\alpha) |= | \frac{(\alpha + \lambda_n a_n)}{(1 + \overline{\lambda_n a_n} \alpha)}|. $
\vskip 0.2in

  Then
 $$
 (1 - |b_n|^2) = 1 - |\frac{\alpha + \lambda_n a_n}{1 + \overline{\lambda_n  a_n} \alpha}|^2
 = \frac{|a_n|^2 (|\alpha|^2 -1) + (1 - |\alpha|^2) }{|1 + \overline{\lambda_n a_n} \alpha|^2}
 \leq \frac{ 2 (1 - |a_n|^2)}{1 - |\alpha|}
 $$

 It follows from our assumption that the sequence $(b_n)$ is a Blaschke sequence.$\:\:\:\:\#\:$\vskip 0.2in

 For part ii) by our assumption $\widetilde{\phi} = \eta \circ \phi\circ \eta_{-1}$ and so assuming
 $(\widetilde{\phi})_n(c_n) = 0 $ we have  $(\widetilde{\phi})_n \circ \eta (d_n) = 0 , $ where $ \eta (d_n) = c_n.$
 Thus $\eta \circ \phi_n(d_n) = 0 .$ By part i) the sequence $(d_n)$ is a Blaschke sequence.$\:\:\:\:\#\:$

\vskip 0.3in
\begin{theorem}
If $ \phi \in \mathbb{A} $ is parabolic and $ \psi \in \mathbb{A}, $ then $ M_{\psi} U_{\phi} $ is not   Crownover.
\end{theorem}

Proof: Suppose $\phi$ is parabolic.Then $\phi \in \mathbb{C}(\phi_{\imath}) $ or  $\mathbb{C}(\phi_{-\imath})$  so by Lemma $3$ and Lemma $4$ ii) , $\phi_n(c_n)=0 \:\:\:\forall n \Rightarrow \{c_n\}$ is a Blaschke sequence. Therefore by Lemma $4$ i) , $\Psi \circ \phi_n (b_n) = 0 \:\:\: \forall n \Rightarrow \{b_n\}$ is a Blaschke sequence. The result now follows from Theorem $2.$

\vskip 0.4in
\section{Isometries of
Codimension Greater than One.}  \vskip 0.2in

 Recall that if S is an isometry on $H^p,( p \neq 2) $ of codimension $d < \infty,$ then  $S = M_{\Psi} C_{\phi}$ for some $\phi \in \mathbb{A}$ and $\Psi $ a d fold Blaschke product. Our key results for codimension $1$ isometries carry over easily to this setting. Thus \vskip 0.1in
 (i) if $ \widetilde{S} \in B (H^p),$ then $$\widetilde{S} \thickapprox S \Leftrightarrow \widetilde{S} = \rho M_{\Psi\circ \eta} U_{\eta_{-1}\circ\phi\circ\eta}$$
  for some $\eta \in \mathbb{A}$ and a $\rho \in T$ which is determined by $\phi$ and $\rho.$ \vskip 0.1in
 (ii) $\bigcap_n^{\infty} S^n H^p = (0)\Leftrightarrow \phi $ is elliptic.

 \vskip 0.2in
 For (ii) , one observes that  the zeros of $(\Psi\circ \phi_n\:\:\:,\:\: n \geq 0)$ can be written as a union of d Blaschke sequences.

 We now consider isometries S of infinite codimension. These can arise in two ways. S could have the form $M_F C_{\Phi}$ where $\Phi$ is inner and $\Phi  \notin  \mathbb{A}.$ The other possibility is that  $S = M_{\Psi} U_{\phi}$ where $\phi  \in \mathbb{A}$ and $\Psi$ is inner and not a finite Blaschke product.
  We focus on this latter case. \vskip 0.2in
   Note that $\bigcap_n^{\infty} S^n H^p = (0)$ if $\phi$ is elliptic.
   \begin{proposition}
   Suppose $\phi \in \mathbb{A}$ and $\phi$ is parabolic or hyperbolic. Depending on this choice of Blaschke product  $\Psi,$ the isometry
   $$
   S = M_{\Psi} U_{\phi}
   $$
   can satisfy either $\bigcap_1^{\infty} S^n H^p = (0)\:$  or   $\:\bigcap_1^{\infty} S^n H^p \neq  (0)$
    \end{proposition}

 \vskip 0.2in
 Proof. The following proof will utilize the results of Theorems $2$ and $4.$ Suppose that $\phi \in \mathbb{A}$ is parabolic and thus by Proposition $1$  we choose $\{\lambda_n\}$ in T   such that $\Psi = \prod_1^{\infty} \lambda_n \phi_{-n}$ is a convergent Blaschke product.

   Note that  $$\Psi \circ \phi = \prod_1^{\infty} \lambda_n \phi_{-n+1} = e \: \Psi. $$
 Iterating this step , we see have $\Psi \circ \phi_n$ has $\Psi$ as a factor $\forall n > 0.$ Thus $\prod_1^{\infty} \Psi \circ \phi_n$ can not be a convergent Blaschke product. Hence,    $ \bigcap_1^{\infty} S^n H^p = \{0\} $ for $S = M_{\Psi} U_{\phi}.$

   Now suppose that $\phi_n(a_n) =0\: \forall n > 0 , $ so that we must have $\sum_1^{\infty} (1 - |a_n|) =R < \infty.$
 Choose $ 1<n_1<n_2<.....,$ such that $ \forall k \geq 1,$
 $$
 \sum_{n=n_k}^{\infty} (1 - |a_n|) < R/(2^k)
 $$
 and let $\Psi = \prod_{k=1}^{\infty} \phi_{n_{k}}.$ Since
 $$
 \sum_{k=1}^{\infty}\sum_{n=n_{k}}^{\infty}(1 - |a_n|) < \sum_{k=1}^{\infty} R/(2^k) < \infty,
 $$

 we see that the zeroes of $\prod_1^{\infty} \Psi \circ \psi_n$ form a Blaschke sequence. Thus $\prod_1^{\infty} \Psi \circ \phi_n$ is convergent, and $\bigcap_1^{\infty} S^n H^p \neq \{0\}$ for $S = M_{\Psi} U_{\phi}.$

\vskip 0.3in
\section{Commuting Automorphisms}\vskip 0.3in

 In this section we elaborate on the conclusion on Corollary $2$ by describing the automorphisms of D that commute
 with a fixed automorphism. These results are undoubtedly known and we outline proofs using the results from
 the last section.
 \begin{definition}
 For $\phi \in \mathbb{A},$ let
 $$
 Com(\phi) = \{ \psi \in \mathbb{A} \:\: : \: \psi \circ \phi = \phi \circ \psi \}.
 $$
 \end{definition}

 Clearly $Com(\phi)$ is a subgroup of $\mathbb{A},$ and $Com(\:e\:)  = \mathbb{A}.$

 \begin{proposition}
 Let $\phi \in  \mathbb{A} \:\:, \phi \neq e.$ \vskip 0.1in
  i)  If $\phi \in \mathbb{E}(a) , $ then $Com(\phi) = \mathbb{E}(a) \bigcup \{ e\}  .$\vskip 0.1in
  ii) If $\phi \in \mathbb{P}(w),$ then $Com(\phi) = \mathbb{P}(w) \bigcup \:\{e\}.$\vskip 0.1in
  iii)If $\phi \in \mathbb{H}(w_1 \: , \:w _2\:),$ then $Com(\phi) = \mathbb{H}(w_1\:\:, w_2\:) \bigcup \{e\}.$\vskip 0.1in
  iv) In each of these cases , $Com(\phi)$ is abelian.
  \end{proposition}

 Proof: Observe that $Com(\eta \circ \phi \circ  \eta_{-1}) = \eta \circ  Com(\phi) \circ \eta_{-1}.$
 So it will suffice to consider $ a = 0\: , \: w=1\:, $ and $ w_1 = -1\:,$ and $ w_2 = 1.$

 For i), if $\phi \in \mathbb{E}(0)= \{ \eta_{\lambda} (z) = \lambda z \::\: \lambda \in T \: , \lambda \neq 1\},$
 then $\mathbb{E}(0) \subset Com(\phi).$  But if $\psi \in Com(\phi),$ then $\psi(0) = \psi(\phi(0)) = \phi(\psi(0)),$
 so $\psi(0) =0.$ Thus $\psi \in \mathbb{E}(0)$ or  $\psi = e.$

 For ii) , if
 $$
 \phi \in \mathbb{P}(1) = \{ \phi_c = \frac{ 1 + c - 2 z}{ 1 + i - z} \:\:: c \in T \:,\: c \neq \pm 1 \},
 $$
 then a computation shows that $\mathbb{P}(1)$ is an abelian set. The same argument as in i) shows that if $\psi \in Com(\phi),$ then $\psi = e $ or $\psi$ has 1 as it unique fixed point. So $Com(\phi) = \mathbb{P}(1) \bigcup (\{e\}).$

 For iii), suppose that $$
 \phi \in \mathbb{H}(-1\:, \: 1 )= \{ \psi_r (z) =  \frac{z-r}{1 + r z} \: : -1 < r <0 , or , \: 0 <r <1\}$$
 Observe that $\mathbb{H}( -1\:, \: 1)$ is an abelian set. Let $\psi \in Com(\phi).$ It follows easily that
 $$
 \{ \psi(1) \: , \psi(-1)\: \} = \{ \pm 1\}.$$
 If $\psi(1) = 1$ and $\psi(-1)= -1,$ then $\psi \in \mathbb{H}(-1\:, \: 1)$ as desired. So suppose $\psi(-1)=1$ and
 $\psi(1)= -1.$  Then $ -\psi \in \mathbb{H}(-1\:, \: 1)$ so $\psi(z) = \frac{r_1 - z}{1 - r_1 z}$ for some $r_1.$ A computation will show that $\psi$ does not commute with automorphisms in $\mathbb{H}( -1\:, \: 1).$ Thus $Com(\phi) = \mathbb{H}(-1\: , \: 1) \bigcup \{e\}.\:\:\:\:\#\:\:$ \vskip 0.2in

References\vskip 0.1in

[$1.$] Cowen,C. and MacCluer, B. : Composition Operators on Spaces of Analytic Functions,
Studies in Advanced Mathematics, CRC Press\vskip 0.1in

[$2.$] Crownover, R. : Commutants of shifts on Banach Spaces, Mich. Math. J. ,Vol. $19,$ Issue $3 (1972), 233-247.$\vskip 0.1in

[$3.$]Duren, P., Theory of $H^p$ Spaces: Academic Press, New York and London, $1970.$\vskip 0.1in

[$4.$] De Leeuw,K., Rudin,W.,  and  Wermer, J.  , The isometries of some function spaces, Proc. Amer. Math. Soc.
$11, (1960), 694- 698.$\vskip 0.1in

[$5.$] Fleming, R. J., and  Jamison, J., Isometries on Banach Spaces: \emph{ function spaces}, Mongraphs and Surveys in
Pure and Applied Mathematics,  $ (2003) \: \: 129$.

Chapman \& Hall\vskip 0.1in

[$6.$]Forelli, F. , The isometries of $H^p,$ Canadian Journal of Mathematics, $16 (1964), 7210 728.$\vskip 0.1in

[$7.$]Robinson, J.: Crownover shift operators., J.of Math. Analysis and Applic. , $02/1988; 130(1): 30-38.$\vskip 0.1in

[$8.$]  Shapiro,J.: Composition operators and classical function theory, Universitat: Tracts in Mathematics,
Springer-Verlag, New York, $1993.$\vskip 0.1in

[$9.$] Sz.Nagy, B. and Foias, C. , Harmonic analysis of operators on Hilbert space, North Holland Publishing Company,
Amsterdam $\centerdot$ London.Inc.

\vskip 1.0in

 Univ. of North Carolina, Chapel Hill, N.C. $27599-3250$\vskip 0.2in

  \end{document}